\documentclass[aap,MSNbibl,nameyear,seceqn,dvips]{arximspdf}

\doi{10.1214/12-AAP870} %kopijuoti is PTS
\volume{23}
\issue{4}
\pubyear{2013}
\firstpage{1291}
\lastpage{1317}

\makeatletter
\newcommand{\R}{\mathbb{R}}
\newcommand{\E}{\mathbb{E}}
\newcommand{\diag}{\operatorname{diag}}
\newcommand{\Cone}{\mathrm{Cone}}
\newtheorem{lemma}{Lemma}
\newproclaim{assumption}{Assumption}
\newproclaim{counterexample}{Counterexample}
\newtheorem{theorem}{Theorem}
\newtheorem{proposition}{Proposition}

\newproclaim{definition}{Definition}

\def\tr{\operatorname{tr}}

\makeatother

\begin{document}
\begin{frontmatter}

\title{Positive recurrence of piecewise Ornstein--Uhlenbeck processes and
common quadratic Lyapunov functions\thanksref{T1}}
\thankstext{T1}{Supported in part by NSF Grant EEC-0926308 and an
Oberwolfach Leibniz fellowship.}
\runtitle{Positive recurrence of piecewise OU processes and CQLFs}

\begin{aug}
\author[A]{\fnms{A. B.} \snm{Dieker}\corref{}\ead[label=e1]{ton.dieker@isye.gatech.edu}}
\and
\author[A]{\fnms{Xuefeng} \snm{Gao}\ead[label=e2]{gxf1240@gatech.edu}}
\runauthor{A. B. Dieker and X. Gao}
\affiliation{Georgia Institute of Technology}
\address[A]{H. Milton Stewart School of Industrial\\
\quad and Systems Engineering\\
Georgia Institute of Technology\\
Atlanta, Georgia 30332\\
USA\\
\printead{e1}\\
\phantom{E-mail: }\printead*{e2}} %adresu isvedimo komanda gale!
\end{aug}

% HISTORY:
\received{\smonth{4} \syear{2011}}
\revised{\smonth{4} \syear{2012}}

% ABSTRACT
%
\begin{abstract}
We study the positive recurrence of piecewise Ornstein--Uhlenbeck (OU)
diffusion processes,
which arise from many-server queueing systems with phase-type service
requirements.
These diffusion processes exhibit different behavior in
two regions of the state space, corresponding to ``overload'' (service
demand exceeds capacity)
and ``underload'' (service capacity exceeds demand).
The two regimes cause standard techniques for proving positive
recurrence to fail.
Using and extending the framework of common quadratic Lyapunov
functions from the theory of control,
we construct Lyapunov functions for the diffusion approximations
corresponding to systems with and without abandonment.
With these Lyapunov functions, we prove that piecewise OU processes
have a unique stationary distribution.
\end{abstract}

% KEYWORDS
\begin{keyword}[class=AMS]
\kwd{60K25}
\kwd{93E15}
\kwd{60J70}
\end{keyword}
\begin{keyword}
\kwd{Stability}
\kwd{common quadratic Lyapunov function}
\kwd{Lyapunov function}
\kwd{piecewise OU process}
\kwd{multi-server queues}
\kwd{customer abandonment}
\kwd{Halfin--Whitt regime}
\kwd{phase-type distribution}.
\end{keyword}

\end{frontmatter}

%s1 #&#
\section{Introduction}
\label{secintro}

Since the pioneering paper of \citet{HalWhi81}, and particularly
within the
last 10 years, there has been a surge of interest in diffusion
approximations for queueing systems
with many servers. These queueing systems model customer contact
centers with hundreds of servers. Empirical study in \citet{BMSZZS05}
suggests that the
service time distribution is far from exponential.
Despite past and foreseeable advances in computer hardware and architectures,
the sheer size of such systems prohibits exact (numerical)
calculations even when the arrival process is Poisson and the service
time distribution is of phase type.
Diffusion approximations such as piecewise Ornstein--Uhlenbeck (OU) processes
can be used to approximate the queue length process.
Such approximations are rooted in many-server heavy traffic limits
proved in
\citet{PuhalskiiReiman00} and \citet{DaiHeTezcan10}. These
approximations are remarkably accurate in predicting system performance
measures, sometimes for systems with as few as $20$ servers [see \citet
{DaiHe11}].

For a diffusion approximation to work, it is critical to know whether
the approximating diffusion process has a unique stationary distribution.
In this paper we prove that, under some natural conditions,
this is the case for piecewise OU processes.
Thus, this paper provides a solid mathematical foundation for \citet{DaiHe11},
who devise an algorithm to numerically compute the stationary
distribution of a piecewise OU process.

A standard technique for
proving stability of queueing systems is to first establish the stability
of a so-called fluid model and then to appeal to general theory
for establishing stochastic stability [see, e.g., \citet{DupuisandWilliams94},
\citet{dai95a},
\citet{sto95}].
However, this theory
is restricted to systems with nonnegative fluid levels which
are attracted to the origin.
The fluid analog of a piecewise Ornstein--Uhlenbeck process does not
possess this property.
As an alternative to the fluid model framework,
the family of \textit{quadratic} Lyapunov functions is a
natural choice for establishing positive recurrence. Indeed,
due to diffusive properties of piecewise Ornstein--Uhlenbeck processes,
if a quadratic Lyapunov function can be shown to stabilize the fluid model,
it simultaneously and directly establishes
stochastic stability, that is, the positive
recurrence of piecewise OU processes.
%Thus, in the context of piecewise OU processes,
%it is `natural' to work with quadratic Lyapunov functions.
As a result of working with quadratic forms as Lyapunov functions,
several key results from linear algebra lie at the heart of
our main results.
We were unable to devise an equally powerful
framework without using this algebraic machinery.

Piecewise OU
processes exhibit different behavior in two regions of the state
space, corresponding to ``overload'' and
``underload.'' The two regions are separated by a hyperplane, which
corresponds to ``critical load.'' In each of the two regions, a~piecewise
OU process can be thought of as a first-order linear
differential equation with stochastic noise. A standard technique in
proving its positive recurrence is to use a quadratic Lyapunov
function to prove stability of such first-order linear differential
equations. However, the two different regions of a piecewise OU
process pose considerable challenges to
apply this methodology. A natural
approach would be to ``paste together'' two quadratic Lyapunov
functions from the two regions, but our attempts in this direction
have failed. In fact, it is well known that a diffusion with two
stable regimes can lead to an instable hybrid system [see
\citet{yinbook} for related examples].
In \citet{blondeltsitsiklis2000},
the stability of a switched linear system is
discussed from the perspective of complexity theory.

%This paper constructs Lyapunov functions for piecewise OU diffusion
%processes
%under assumptions that are natural in view of
Using the interpretation of the diffusion parameters in terms of a
many-server queueing system, our main results can be formulated as
follows: (1) For a slightly underloaded system without abandonment,
we show that there exists a quadratic Lyapunov function which yields
the desired positive recurrence using the Foster--Lyapunov criterion
(Theorem~\ref{thmnoabd}). In general, this quadratic Lyapunov
function is not explicit and nonunique. (2) We show that no
quadratic Lyapunov function can satisfy the Foster--Lyapunov
criterion for systems with abandonment. (3) We construct a suitable
nonquadratic Lyapunov function to prove positive recurrence for
systems with abandonment (Theorem~\ref{thmqedabd}).

The main building blocks for these two types of Lyapunov functions
are so-called common quadratic Lyapunov functions (CQLFs), which are
widely used in the
theory of control. Such functions play an
important role in the stability analysis for deterministic
linear systems, with different dynamics in different parts of the state space
(or, more generally, operating under a switching rule).
They are called \textit{common} quadratic Lyapunov functions
since they serve as a quadratic Lyapunov function
in each part of the state space. There is a vast body of literature on
CQLFs and related
theory [see the survey \citet{Shorten07} for details]. Although
quadratic Lyapunov functions are ubiquitous in the literature on
queueing systems \citet{daipra00},
\citet{GamarnikMomcilovic08},
\citet{taseph92}, to
our knowledge, our paper is the first to exploit CQLFs in this
context.
%gamarnik goldberg?

As mentioned in the section on open problems of \citet{Shorten07},
it is of considerable interest to determine simple conditions for the
existence of CQLFs.
Theorem~\ref{thmexistcqlf}, which is our main technical contribution
in this space,
establishes such a result in the context of M-matrices and rank-1 perturbations.
The theorem shows that existence of a CQLF is guaranteed after merely
verifying that certain vectors
are nonnegative. It is a first result of this kind.
Its proof relies on a delicate analysis involving Chebyshev polynomials,
as well as on an extension of recent work of \citet{King06} and \citet
{Shorten09}
summarized in Proposition~\ref{propshorten2009} below.

To conclude this \hyperref[secintro]{Introduction}, we mention
a body of work on the recurrence of multidimensional
Ornstein--Uhlenbeck type processes by \citet{Sato94},
\citet{Sato96}, which differ from the processes
studied here.
%through Fourier-analytic methods.
The processes studied in these references are driven by L\'evy
processes and they have a linear
drift coefficient. %, which are associated with diagonalizable matrices
%or Jordan cell matrices.
As a result, their multidimensional processes
do not possess the critical feature of the processes we study here,
namely, a piecewise linear drift coefficient.

This paper is organized as follows.
Section~\ref{sec2} discusses the required background on piecewise OU process and
positive recurrence.
Section~\ref{secCQLF} is devoted to common quadratic Lyapunov functions.
Section~\ref{secmainre} summarizes the main results and Section~\ref{sec5} contains the proofs
of the main results.
The proof of Proposition~\ref{propshorten2009}, which mainly uses
existing methodology from the theory of control,
is given in Appendix~\ref{appendixprop3pf}.
Appendix~\ref{appendixcountereg} shows that no quadratic Lyapunov function can work in the
Foster--Lyapunov criterion
if abandonment is allowed.

\subsection*{Notation}
All random variables and stochastic processes are defined on a common
probability space $(\Omega,\mathcal{F},\mathbb{P})$ unless otherwise
specified.
%The symbols $\mathbb{Z}$, $\mathbb{Z}_{+}$, $\mathbb{N}$, $
%are used to denote the sets of integers, nonnegative integers,
%positive integers, real numbers, and nonnegative real numbers,
%respectively.
For some $d\in\mathbb{N}$, $\mathbb{R}^{d}$ denotes the
$d$-dimensional Euclidean space. %; thus, $\mathbb{R}=\mathbb{R}^{1}$.
The space of functions $f\dvtx\mathbb{R}^{K}\rightarrow\mathbb{R}$ that
are twice continuously differentiable is denoted by
$C^2(\mathbb{R}^{K})$. We use $\nabla$ to denote the gradient
operator. %and $f^{''}$ to denote the second derivative, which is a $K
Given $x\in\mathbb{R}$, we set $x^{+}=\max\{x,0\}$.
%$x^{-}=\max\{-x,0\}$, and $\lfloor x\rfloor=\max\{j\in\mathbb{Z}:j\leq
%x\}$.
All vectors are envisioned as column vectors. For a $K$-dimensional
vector $u$, we use $u_{k}$ to
denote its $k$th entry and we write $\vert u\vert$ for its Euclidean
norm. We also
write $u'$ for its transpose. For two $K$-dimensional vectors $u$
and $v$, we write $u' \ge v'$ $(u'>v')$ if $u_{k}\ge v_{k}$
$(u_{k}>v_{k})$ for each $k=1,2,\ldots, K.$ The inner product of $u$
and $v$ is denoted by $u'v,$ which is $\sum_{k = 1}^K u_k
v_k.$ Given a $K \times K$ matrix $M$, we use $M^{\prime}$ to denote
its transpose and $M_{ij}$ for its $(i,j)$th entry. We write $M>0$
$(M<0)$ if $M$ is a positive (negative) definite matrix and $M \ge
0$ $(M \le0)$ if it is a positive (negative) semi-definite matrix.
Let the matrix norm of $M$ be $|M|=\sum_{ij}
|M_{i,j}|,$ where $|M_{ij}|$ is the absolute value of
$M_{ij}.$ We reserve $I$ for the $K\times K$ identity matrix and $e$
for the $K$-dimensional vector of ones.

%s2 #&#
\section{Piecewise OU processes and positive recurrence}\label{sec2}
This section introduces the piecewise Ornstein--Uhlenbeck (OU)
processes studied in this paper, and discusses preliminaries on
positive recurrence.

%s2.1 #&#
\subsection{Piecewise OU processes} \label{subsecpiecewiseou}

We first define M-matrices. We call a matrix nonnegative when each
element of the matrix is nonnegative.

%de1 #&#
%
\begin{definition}[(M-matrix)]
A $K \times K$ matrix $R$ is said to be an M-matrix if it can be
expressed as $R=sI-N$ for some $s>0$ and some nonnegative matrix~$N$
with the property that $\rho(N) \le s$, where $\rho(N)$ is the
spectral radius of $N.$ The matrix $R$ is nonsingular if $\rho(N) <
s$.
\end{definition}

We next define piecewise Ornstein--Uhlenbeck (OU) processes, which
are special diffusion processes. Let $\{W(t)\}$ be a standard Brownian
motion in any dimension. A $K$-dimensional diffusion process $Y$ is
the strong solution to a stochastic differential equation of the
form
\[
dY(t)=b(Y(t))\,dt+ \sigma(Y(t)) \,dW(t),
\]
where the drift coefficient $b(\cdot)$ and the diffusion coefficient
$\sigma(\cdot)$ have appropriate sizes and satisfy the following
Lipschitz continuity condition: there exists some $C>0$ such that
%e1 #&#
%
\begin{equation} \label{eqdifcoeflip}
|b(x) - b(y)| + |\sigma(x) - \sigma(y)| \le C|x - y|\qquad
\mbox{for all $x, y \in \mathbb{R}^{K}.$ }
\end{equation}
For a real-valued function $V \in C^2(\mathbb{R}^{K})$, the
generator $G$ of $Y$ applied to $V$ is given by, for $y \in\mathbb
{R}^{K},$
%e2 #&#
%
\begin{equation} \label{eqgenerator}
GV(y) = (\nabla V (y) )'b(y) + \frac{1}{2}\sum_{i,j}
{(\sigma\sigma')_{ij} (y)\frac{{\partial^2 V}}{{\partial y_i\,
\partial y_j }}}(y).
\end{equation}
We refer to \citeauthor{RogersWilliamsvol2} [(\citeyear{RogersWilliamsvol2}), Chapter V], for more details
on diffusion processes.

%de2 #&#
%
\begin{definition}[(Piecewise OU processes)] \label{defpiecewiseou}
Let $p$ be a $K$-dimensional probability vector, $e$
be the $K$-dimensional vector of ones and let $R$ be a $K \times
K$ nonsingular M-matrix. For $\alpha, \beta\in\mathbb{R}$, a
$K$-dimensional diffusion process $Y$ is called a piecewise
Ornstein--Uhlenbeck $($OU$)$ process if it has drift coefficient
%e3 #&#
%
\begin{equation}
\label{eqdrfitcoef} b(y)= - \beta p -R\bigl(y-p(e'y)^+ \bigr) -\alpha p
(e'y)^+ ,
\end{equation}
and diffusion coefficient $\sigma(y) \equiv\sigma$ for all $y \in
\mathbb{R}^K$, such that $\sigma\sigma'$ is a $K \times K$
nonsingular matrix.
\end{definition}

As in \citet{DaiHeTezcan10}, we call this process a \textit
{piecewise} OU process
since the drift coefficient is affine (hence, OU process) yet
it differs on each side of the hyperplane $\{y\in\R^K\dvtx e'y=0\}$
(hence, piecewise).
Indeed, for $e'y\ge0$ we have $b(y)=-\beta p -R(I-pe')y - \alpha p (e'y)$
while for $e'y\le0$ we have $b(y) = -\beta p -Ry$.
In conjunction with $\sigma(y) \equiv
\sigma$, this implies the Lipschitz continuity condition
(\ref{eqdifcoeflip}). As a consequence, the piecewise OU process
$Y$ is well-defined as a diffusion process.

The quantities $\alpha,\beta, R,p$ on the right-hand side of
(\ref{eqdrfitcoef}) come from the queueing system that gave rise to
the piecewise OU diffusion. Their queueing interpretation is as
follows: $\alpha$ is the abandonment rate, $\beta$ is the slack in
the arrival rate relative to a critically loaded system while $p$
and $R$ are the parameters of the service-time distribution (assumed
to be of phase-type). For more details, we refer to
\citet{DaiHeTezcan10}.

Throughout the paper, we impose the following assumption.

%as1 #&#
%
\begin{assumption}\label{asseR}
Each component of the row vector $e'R$ is nonnegative, that is,
\[
e'R \ge0'.
\]
\end{assumption}

We now make the connection between piecewise
OU processes and many-server queueing models explicit, and we discuss
Assumption~\ref{asseR}
in this context.
For presentational convenience, we do so in the special case of
$M/H_2/n+M$ queues.
In fact, we consider a sequence of
$M/H_2/n+M$ queues indexed by $n$, where $n$ is the number of
(identical) servers, meaning that (1) the arrival process is
a Poisson process with some intensity $\lambda^n$, (2) the service times
have a two-phase hyperexponential distribution, so they are
exponential with parameter $\nu_1$ with probability $p_1$,
and exponential with parameter $\nu_2$ with probability $p_2=1-p_1$ and
(3) each customer has a patience time which follows an exponential
distribution with parameter $\alpha>0$.
Hypergeometric service time distributions are of special interest,
since they can be used to model
multiclass systems\vadjust{\goodbreak} [see \citet{PuhalskiiReiman00},
\citet{gamarnikstolyar2011}].
To see this, envision
two types of customers entering a buffer to seek service. Suppose that
the mean
service time is 1, that is, ${p_1/\nu_1 + p_2/ \nu_2 }=1$. We further
assume the system is operated under Halfin--Whitt regime, that is, for
some $\beta\in\R$,
\[
\lim_{n \rightarrow\infty}\sqrt{n} \biggl(1- \frac{\lambda^n}{n}\biggr) =
\beta.
\]
Let $X_1^n(t)$ and $X_2^n(t)$ denote the number of customers
of type 1 and 2 in the system at time $t$. For $i=1, 2$ and $t
\ge0,$ we define
\[
\tilde X_i^n (t) = \frac{1}{\sqrt{n}} \biggl(X_i^n(t)- n\frac{p_i}{ \nu
_i} t\biggr).
\]
As detailed in \citet{DaiHeTezcan10} [see \citet{GamarnikGoldberg11}
for a related general result],
the ``centering'' in this expression has been chosen so that,
in a sense of weak convergence on the process level,
\[
( \tilde X_1^n, \tilde X_2 ^n) \Rightarrow(Y_1, Y_2),\qquad n \rightarrow
\infty,
\]
where $(Y_1, Y_2)$ satisfies the following system of stochastic
differential equations:
for $i=1, 2,$
\begin{eqnarray*}
Y_i(t) &=& Y_i(0) +W_i(t)- \beta p_i t \\
&&{}- \nu_i \int_{0}^{t}
\bigl(Y_i(s)- p_i \bigl(Y_1(s)+ Y_2(s)\bigr)^+\bigr) \,ds - \alpha p_i
\int_{0}^{t}\bigl(Y_1(s)+ Y_2(s)\bigr)^+ \,ds.
\end{eqnarray*}
Note that $(Y_1(s)+Y_2(s))^+$ represents the (scaled) number of
customers waiting in the buffer,
and the fraction of type $i$ customers in the buffer is approximately
$p_i$. Thus,
the term involving $\nu_i$ can be thought of as a service-rate term.
Similarly, the terms involving $\alpha$ and $\beta$ are the
abandonment and arrival term,
respectively.
The randomness in the system is represented by $W=(W_1, W_2)$, which is a
driftless Brownian motion with nonsingular covariance matrix
\[
\pmatrix{
{p_1(p_1 c^2 - p_1 +2)} & {p_1 p_2 (c^2 -1)} \vspace*{2pt}\cr
{p_1 p_2 (c^2 -1) } & {p_2 (p_2 c^2 -p_2 +2)} }
\]
for some constant $c\in\R$. Therefore, $Y=(Y_1, Y_2)$ is a
two-dimensional piecewise OU process with drift coefficient
%e4 #&#
%
\begin{equation}
b(y)= - \beta p - R \bigl(y-p(e'y)^+ \bigr) -\alpha p (e'y)^+ ,
\end{equation}
where the matrix $R$ is given by
\[
R=\pmatrix{
\nu_1 & 0 \vspace*{2pt}\cr
0 & \nu_2 }.
\]

When we apply this procedure to a general
phase-type service time distribution with $K$ phases, the
corresponding diffusion limit is a $K$-dimensional piecewise OU
process $Y$. The parameters $p$ and $R$ represent the distribution
of the initial phase and phase dynamics, respectively. Each
component of the piecewise OU process $Y_k$ approximates the number
of phase-$k$ customers in the many-server queueing system, either
waiting or in service. Thus $e'Y$ represents the total number of
customers in the system after centering and scaling. Thus, whenever
$e'Y>0$ the system is in ``overload,'' that is, there are customers
waiting in the buffer, and whenever $e'Y<0$ the system is in
``underload,'' that is, there are idle servers. We refer readers
to~\citet{PuhalskiiReiman00} and \citet{DaiHeTezcan10} for more
details. We remark that the matrix $R$ in these two papers takes the
form of $(I-P') \diag\{\nu\},$ where $P$ is assumed to be a
transient matrix describing the transitions between each service
phase, and $\diag\{\nu\}$ is a diagonal matrix with $k$th
diagonal entry given by $\nu_k$, where $\nu_k$ is the rate for the
sojourn time in phase $k$. Transience of $P$ corresponds to
customers who eventually leave the system after receiving a
sufficient amount of service, which implies that $e'R =e'(I-P')\diag
\{\nu\} \ge0$. Therefore, we conclude that in this setting, $R$ is
a nonsingular M-matrix and that Assumption 1 is satisfied.

%s2.2 #&#
\subsection{Positive recurrence and Lyapunov functions}
\label{secposrec} In this section, we recall the definitions and
the criteria for positive recurrence and exponential ergodicity in
the context of general diffusion processes.

Let $\E_\pi$ be the expectation operator with respect to a
probability distribution $\pi$.

%de3 #&#
%
\begin{definition}[(Positive recurrence and stationary distribution)]
For a $K$-dimensional diffusion process $Y,$ we say that $Y$ is
positive recurrent if for any $y \in\mathbb{R}^{K}$ and any
compact set $C$ in $\mathbb{R}^{K}$ with positive
Lebesgue measure, we have
\[
\E\bigl(\tau_C | Y(0)=y\bigr)< \infty,
\]
where $\tau_C = \inf\{t \ge0\dvtx Y(t) \in C\}$ is the hitting time of
the set $C$. We call a probability distribution $\pi$ on $\mathbb
{R}^{K}$ a stationary distribution for $Y$ if for every bounded
continuous function $f$: $\mathbb{R}^{K} \rightarrow\mathbb{R}$,
\[
\mathbb{E}_\pi[f(Y(t))]=\mathbb{E}_\pi[f(Y(0))]\qquad \mbox{for
all $t \ge0.$}
\]

\end{definition}

In the following, we assume that the diffusion coefficient of the
diffusion process $Y$ is uniformly nonsingular. That is, there
exists some $c \in(0, \infty)$ such that for all $y \in
\mathbb{R}^{K}$ and $a \in\mathbb{R}^{K}$,
%e5 #&#
%
\begin{equation}\label{eqNDcoef}
a' \sigma(y) \sigma(y)' a \ge c a'a.
\end{equation}

The next result gives a sufficient criterion for
positive recurrence of diffusion processes [see \citet
{Khasminskii2011}, Sections 3.7, 4.3 and 4.4 and \citet{meytwe93b}, Section~4].
Uniqueness of the stationary distribution follows from
\citet{Peszat95} in view of condition (\ref{eqNDcoef}).

%pr1 #&#
%
\begin{proposition}[(Foster--Lyapunov criterion)]
\label{propprec}
Let $Y$ be a diffusion process satisfying
$(\ref{eqNDcoef})$. Suppose that there exists a nonnegative function
$V \in C^2(\mathbb{R}^{K})$ and some $r>0$ such that, for any $|y|
> r,$
\[
GV(y) \le- 1.
\]
In addition, suppose that $V(y) \to\infty$ as $ |y| \to\infty$.
Then $Y$ is positive recurrent and has a unique stationary
distribution. The function $V$ is called a Lyapunov function.
\end{proposition}

We now introduce the concept of exponential ergodicity. For any
positive measurable function $f \ge1$ and any signed measure $m$,
we write $\Vert m\Vert_f=\sup_{|g|\le f} |m(g)|$.

%de4 #&#
%
\begin{definition}[(Exponential ergodicity)]
Suppose that the diffusion process~$Y$ is positive recurrent and
that it has a unique stationary distribution~$\pi$. Given a function
$f \ge1$, we say that~$Y$ is $f$-exponentially ergodic if there
exists a $\gamma\in(0,1)$ and a real-valued function $B$ such that
for all $t>0$ and $y\in\mathbb{R}^K$,
\[
\Vert P^t(y, \cdot)- \pi(\cdot)\Vert_f \le B(y) \gamma^t,
\]
where $P^t$ is the transition function of $Y.$ If $f\equiv1$, we
simply say that $Y$ is exponentially ergodic.
\end{definition}

For $f \ge1,$ we have $\Vert P^t(y, \cdot)- \pi(\cdot)\Vert_1 \le
\Vert P^t(y,
\cdot)- \pi(\cdot)\Vert_f,$ and we deduce that $f$-exponential
ergodicity implies exponential ergodicity. The following result
gives a criterion for exponential ergodicity [see
\citet{meytwe93b}, Section~6].

%pr2 #&#
%
\begin{proposition}
\label{propexpergo} Suppose that $Y$ is a diffusion process with a
unique stationary distribution. If there is a nonnegative function
$V \in C^2(\mathbb{R}^{K})$ such that $V(y) \to\infty$ as $ |y|
\to\infty$ and for some $c>0, d< \infty$,
\[
GV(y) \le-c V(y)+d \qquad\mbox{for any $y \in\mathbb{R}^K$,}
\]
then $Y$ is $(V+1)$-exponentially ergodic.
\end{proposition}

%s3 #&#
\section{Common quadratic Lyapunov functions}
\label{secCQLF} In this section we introduce common quadratic
Lyapunov functions (CQLFs). Such functions play a central role in
the stability analysis of deterministic switched linear systems,
which is discussed in Section~{\ref{seccqlfexist}}. We use CQLFs as
building blocks to construct Lyapunov functions to prove positive
recurrence of piecewise OU processes. At this point it is best to
distinguish CQLFs for switched linear systems from the Lyapunov
functions in the context of the Foster--Lyapunov criterion. We
connect these two concepts in Section~\ref{secmainre}.

%s3.1 #&#
\subsection{Background and definitions}\label{seccqlfdef}
Quadratic Lyapunov functions form a cornerstone of stability theory
for ordinary differential equations. Consider the linear system
$\dot y(t)=B y(t)$ where $y(t) \in\mathbb{R}^{K}$, $B \in
\mathbb{R}^{K\times K}$ is a fixed real matrix and $\dot y(t)$ is
the derivative of $y$ with respect to $t$. For $Q\in
\mathbb{R}^{K\times K}$, the quadratic form $L$ given by $L(y)=y'Qy$
for $y \in\mathbb{R}^K$ is called a quadratic Lyapunov function for
the matrix~$B$ if $Q$ is positive definite and $QB+B'Q$ is negative
definite. In this case, there exists a constant $C>0$ such that
\[
\frac{d}{dt}{L(y(t))}=y(t)' (QB+B'Q) y(t)\le-C L(y(t)) <0\qquad
\mbox{for all $t \ge0,$}
\]
and thus we can conclude that
$L(y(t)) \le e^{-Ct} L(y(0)).$ This implies that $L(y(t))\rightarrow
0$ as $t \rightarrow\infty,$ thus $y(t)\rightarrow0$ as $t
\rightarrow\infty.$ It is standard fact in Lyapunov stability
theory that the existence of a quadratic Lyapunov function $L$ is
equivalent to all eigenvalues of $B$ having negative real part [\citet
{bermanplemmons}, Section~6.2].

The following definition, tailored to our setting in order to allow
for a singular matrix, plays an important role in our analysis.
Other versions can be found in \citet{ShortenNarendra03} and
\citet{Shorten07}. Recall that an eigenvalue of a matrix is called
(geometrically) simple if its corresponding eigenspace is
one-dimensional.

%de5 #&#
%
\begin{definition}[(CQLF)] \label{defcqlf}
Let $B_1 \in\mathbb{R}^{K\times K}$ have all eigenvalues with
negative real part and let $B_2 \in\mathbb{R}^{K\times K}$ have all
eigenvalues with negative real part except for a simple zero
eigenvalue. For $Q\in\mathbb{R}^{K\times K},$ the quadratic form
$L$ given by $L(y)=y'Qy$ for $y \in\mathbb{R}^K$ is called a common
quadratic Lyapunov function $($CQLF$)$ for the pair $(B_1, B_2)$ if $Q$
is positive definite and
\begin{eqnarray*}
QB_1+B_1'Q&<&0, \\
QB_2+ B_2'Q &\le&0.
\end{eqnarray*}
\end{definition}

%s3.2 #&#
\subsection{The CQLF existence problem} \label{seccqlfexist}
The CQLF existence problem for a pair of matrices has its roots in
the study of stability criteria for switched linear systems. These
systems have the form $\dot y(t)=B (\tau)y(t)$, where $B(\tau)\in
\{B_1, B_2\} $ with $B_i \in\mathbb{R}^{K\times K}$ for $i=1, 2$
and where the switching function $\tau$ may depend on both $y$ and
$t.$ The existence of a CQLF for the pair $(B_1, B_2)$ guarantees
that all solutions of the systems are bounded under arbitrary
switching function $\tau.$ The CQLF existence problem is also
closely related to the Kalman--Yacubovich--Popov lemma in the
development of adaptive control algorithms and the Lur'e problem in
nonlinear feedback analysis. For more details consult
\citet{KYPlemma}, \citet{LMILur} and the recent survey paper by
\citet{Shorten07}. For an arbitrary matrix pair, no simple analytic
and verifiable conditions are known for the pair to admit a CQLF.
In the special case where the difference of the
matrices has rank one, \citet{King06}\vadjust{\goodbreak} shows that if both $B_1$ and
$B_2$ are Hurwitz, that is, all eigenvalues of the matrices $B_1, B_2$
have negative real part, then there exists a positive definite
matrix $Q$ such that $QB_1+B_1'Q<0$ and $QB_2 +B_2'Q<0$ if and only
if the matrix product $B_1 B_2$ has no real negative eigenvalues.
Note that in this case, both $B_1$ and $B_2$ are nonsingular. A similar
CQLF existence result has been obtained in \citet{Shorten09} when
one of the matrices ($B_1$ or $B_2$) is singular.

We now state a result on the CQLF existence problem
for a pair of matrices with one of them being singular. It is
essentially the main theorem in \citet{Shorten09} but we relax their
assumptions. Let $B \in\mathbb{R}^{K\times K}$ be a real matrix
and let $g, h \in\mathbb{R}^K$. The proposition below is stated in
\citet{Shorten09} under the assumptions that $(B,g)$ is controllable,
meaning that the vectors $g,Bg,B^2 g, \ldots$ span $\mathbb{R}^K$,
and that $(B,h)$ is observable, meaning that the vectors
$h,B'h,(B')^2 h, \ldots$ span $\mathbb{R}^K$. Using techniques from
\citet{King06}, we show that these assumptions are unnecessary and we
state the result in its full generality here. A~proof is given in
Appendix~\ref{appendixprop3pf}.

%pr3 #&#
%
\begin{proposition}
\label{propshorten2009} Suppose that all eigenvalues of matrix $B$
have negative real part and all eigenvalues of $B-gh'$ have negative
real part, except for a simple zero eigenvalue. Then there exists a
CQLF for the pair $(B, B-gh')$ if and only if the matrix product $B
(B-gh')$ has no real negative eigenvalues and a simple zero
eigenvalue.
\end{proposition}

%s4 #&#
\section{Main results} \label{secmainre}
In this section, we present our results on positive recurrence of
the piecewise OU process $Y$. Key to these results is the following
theorem, which uses Proposition~{\ref{propshorten2009}} to
establish the existence of a CQLF for certain matrix pairs. Recall
the definitions of $R,$ $p$ and $e$ from
Definition~\ref{defpiecewiseou} in
Section~\ref{subsecpiecewiseou}, and note that we are working under
Assumption~\ref{asseR}.

%th1 #&#
%
\begin{theorem}
\label{thmexistcqlf} There exists a CQLF for both the pair $(-R,
-R(I-pe'))$ and the pair $(-R, -(I-pe')R).$
\end{theorem}

By Theorem~\ref{thmexistcqlf}, there exists a CQLF $L$ for the pair
$(-R, -R(I-pe'))$ and another CQLF $\widetilde L$ for the pair $(-R,
-(I-pe')R).$ Typically there are many CQLFs corresponding to these
pairs, that is, $L$ and $\widetilde L$ are not unique. Note that $L$
and $\widetilde L$ are closely related in the following sense. If
the CQLF $L$ for the pair $(-R, -R(I-pe'))$ is given by $L(y)=y'Qy$
for some $Q>0$ and for all $y \in\mathbb{R}^K,$ then one readily
checks that the quadratic form $\widetilde L$ given by $\widetilde
L(y)=y'(R'QR)y $
for $y \in\mathbb{R}^K$ is a CQLF for the
pair $(-R, -(I-pe')R).$ We remark that, apart from special cases,
the CQLFs from Theorem~\ref{thmexistcqlf} are not
explicit.\looseness=-1

We know from Theorem~\ref{thmexistcqlf} that there exists a CQLF
$L$ for the pair $(-R, -R(I-pe')),$ where $L$ is given by
$L(y)=y'Qy$ for some $Q>0$ and for all $y \in\mathbb{R}^K.$ We are
able to use the quadratic form $L$ as a Lyapunov function in the
Foster--Lyapunov criterion of Proposition~\ref{propprec} to prove
the following result.\vadjust{\goodbreak}

%th2 #&#
%
\begin{theorem}
\label{thmnoabd} If $\alpha= 0 $ and $ \beta
> 0$, then the piecewise OU process $Y$ is positive
recurrent and has a unique stationary distribution.
\end{theorem}

For $\alpha>0$, no quadratic function can serve as a Lyapunov
function in the Foster--Lyapunov criterion to prove positive
recurrence of the piecewise OU process $Y$ (see Appendix
\ref{appendixcountereg} for details). Despite this fact, still
relying on Theorem~\ref{thmexistcqlf}, we overcome this difficulty
in Section~\ref{secpfofthm3}\vspace*{1pt} by constructing a suitable
nonquadratic Lyapunov function. Specifically, there exists a CQLF
$\widetilde L$ for the pair $(-R, -(I-pe')R)$ by
Theorem~\ref{thmexistcqlf}, where $\widetilde L$ is given by
$\widetilde L(y)=y'\widetilde Qy$ for some $\widetilde Q>0$ and for
all $y \in\mathbb{R}^K.$ A suitable approximation to the function
$f$, given by, for all $y\in\mathbb{R}^K,$
\[
f(y)=(e'y)^2+ \kappa\widetilde L \bigl(y- p (e'y)^+\bigr)\qquad \mbox{for
some large constant $\kappa$,}
\]
provides the desired nonquadratic Lyapunov function in the
Foster--Lyapunov criterion to prove positive recurrence of $Y$ when
$\alpha>0$.
Note that, in queueing terminology, the vector $y-p(e'y)^+$ relates to the
customers in service, and not to those in the buffer.
We therefore need the extra term $(e'y)^2$.
Applying Proposition~\ref{propexpergo} with
the same nonquadratic Lyapunov function yields exponential
ergodicity of $Y$ for $\alpha>0.$ We use a smooth approximation of
$f$ as a Lyapunov function in the Foster--Lyapunov criterion of
Proposition~\ref{propprec} instead of using $f$ directly since $f
\notin C^2(\mathbb{R}^K)$. This leads to the following result.

%th3 #&#
%
\begin{theorem}
\label{thmqedabd} If $\alpha>0$, then the piecewise OU process $Y$
is positive recurrent and has a unique stationary distribution.
Moreover, $Y$ is exponentially ergodic.
\end{theorem}

%s5 #&#
\section{Proof of the main results}\label{sec5}

%s5.1 #&#
\subsection{\texorpdfstring{Proof of Theorem~\protect\ref{thmexistcqlf}}{Proof of Theorem 1}}
\mbox{}
\begin{pf}We only establish the existence of a CQLF for the pair
$(-R,\break
-R(I-pe')),$ since the existence of a CQLF for the other pair
$(-R,
-(I-pe')R)$ follows directly. Since $-R-
(-R(I-pe')) =-R pe'$ is a rank-one matrix, in view of
Proposition~\ref{propshorten2009}, we need to check three
conditions:
\begin{longlist}[(a)]
\item[(a)] All eigenvalues of $-R$ have negative real part.
\item[(b)] All eigenvalues of $-R(I-pe')$ have negative real part except
for a simple zero eigenvalue.
\item[(c)] The matrix product $R^2(I-pe')$ has no real negative
eigenvalues and a simple zero eigenvalue.
\end{longlist}

We first prove (a) and (b). It is known that all eigenvalues of
a nonsingular M-matrix have positive real part, and all eigenvalues
of a singular M-matrix have\vadjust{\goodbreak} nonnegative real part [see
\citet{bermanplemmons}, Chapter~6]. Since $R$ is a nonsingular
M-matrix, we immediately get (a). For (b), it is clear that
$-R(I-pe')$ has a simple zero eigenvalue. We notice that
$(I-pe')R=R-p e'R$ where $e'R \ge0'$ by Assumption~\ref{asseR},
$p$ is a nonnegative vector and $R$ is a nonsingular
M-matrix, so the off-diagonal elements of $(I-pe')R$ are
nonpositive. Using this in conjunction with the fact that both
$I-pe'$ and $R$ are M-matrices, we find that $(I-pe')R$ is also an
M-matrix and all its eigenvalues have nonnegative real part [see
\citet{bermanplemmons}, Exercise~5.2]. Thus we get (b) after a
similarity transform.

We now concentrate on proving (c). The key ingredient of the proof
is an identity for Chebyshev polynomials. Suppose that $R^2(I-pe')$
has a real negative eigenvalue $-\lambda$ with $\lambda>0$, and
write $v$ for the corresponding left eigenvector, thus we have
$v'R^2(I-pe')=-\lambda v'$. Right-multiplying by $p$ on both sides,
we obtain $v'p=0$ and the following equality:
%e6 #&#
%
\begin{eqnarray}\label{eqveigen}
0&=&v'R^2(I-pe')+\lambda v'=v'R^2(I-pe')+\lambda
v'(I-pe')
\nonumber
\\[-8pt]
\\[-8pt]
\nonumber
&=&v'(R^2+\lambda I)(I-pe').
\end{eqnarray}
Since $R$ is a nonsingular M-matrix having only eigenvalues with
positive real part, the matrix $(R^2+\lambda I)$ is invertible for
all $\lambda>0$. Also, by the fact that $p$ is a
nonnegative probability vector with $e'p=1,$ we deduce the matrix
$(I-pe')$ has an eigenvalue 0 and the corresponding left eigenvector
must be in the form of $c e'$ for some $c \ne0$. Thus, it follows
from (\ref{eqveigen}) that $v'=c e'(R^2+\lambda I)^{-1}$ for some
$c \ne0$. We show below that $e'(R^2+\lambda I)^{-1}$ is a
positive vector for all $\lambda>0$, that is,
%e7 #&#
%
\begin{equation}\label{eqpositivevec}
e'(R^2+\lambda I)^{-1}>0'\qquad \mbox{for all $\lambda>0$.}
\end{equation}
This yields a contradiction in view of $v'p=0$. By definition of a
nonsingular M-matrix, $R$ is of the form $sI-N$,
where $N$ is a nonnegative matrix with $\rho(N) < s$ and $e'R \ge0$
by Assumption~\ref{asseR}. Inequality~(\ref{eqpositivevec}) thus
states that for all $\lambda>0$ and for every nonnegative matrix $N$
with $\rho(N) < s$ and $se' \ge e'N,$
\[
e'\bigl((sI-N)^2+\lambda I\bigr)^{-1}>0' .
\]
Equivalently, we show the following inequality: for all $y \in(0,
1)$ and for every nonnegative matrix $N$ with $\rho(N) < 1$ and $e'
\ge e'N,$
%e8 #&#
%
\begin{equation}\label{eqpositivecheb}
e'\bigl(y(I-N)^2+ (1-y) I\bigr)^{-1}>0'.
\end{equation}
Therefore, to show (c), it suffices to prove
(\ref{eqpositivecheb}) for fixed $N$ and $y \in(0, 1)$.

Our strategy to prove (\ref{eqpositivecheb}) is to use a matrix
series expansion and connections with Chebyshev polynomials.
Chebyshev polynomials of the second kind $U_n$ can be defined by the
following trigonometric form:
%e9 #&#
%
\begin{equation}\label{eqUtri}
U_n(\cos\theta)=\frac{\sin(n+1) \theta}{ \sin\theta}
\qquad\mbox{for $n=0, 1, 2, 3, \ldots.$}
\end{equation}
Moreover, for $z \in[-1,1]$ and $t \in(-1, 1)$, the generating
function of $U_n$ is
%e10 #&#
%
\begin{equation}\label{eqUgene}
\sum_{n=0}^\infty{ U_n(z) t^n }=\frac{1}{1-2tz+t^2}.
\end{equation}
Refer to \citet{AbramowitzandStegun}, Chapter~22, for more details.
The scalar version of the left-hand side of (\ref{eqpositivecheb})
admits the following expansion: for $ x, y \in(0,1)$,
%e11 #&#
%
\begin{equation}
\label{eqcheb} \frac{1}{{ y{{(1 - x)}^2+ 1 - y}}} = \sum_{n =
0}^\infty{{C_n}(y){x^n}},
\end{equation}
where ${C_n}(y) = {U_n}(\sqrt y ){(\sqrt y )^n}$ for all $n \ge0$.
This can readily be verified with (\ref{eqUgene}). In particular,
we have
%e12 #&#
%
\begin{equation}\label{eqC0}
C_0(y)= U_0(y) \equiv1\qquad \mbox{for all $y \in(0, 1)$.}
\end{equation}
For fixed $y \in(0,1),$ the radius of convergence of the power
series in (\ref{eqcheb}) is larger than 1. Since $\rho(N) < 1$, we
immediately obtain that, for $y \in(0,1),$
%e13 #&#
%
\begin{equation}
\label{eqchebmatrix} \bigl(y(I-N)^2+ (1-y) I\bigr)^{-1} = \sum_{n =
0}^\infty{{C_n}(y){N^n}}.
\end{equation}
Let $y \in(0,1)$ be fixed and define $\theta$ through $\sqrt y
=\cos\theta\in(0, 1).$ Using the trigonometric form
(\ref{eqUtri}) of $U_n$, we can then show by induction that, for any
$m \ge1$,
%e14 #&#
%
\begin{eqnarray} \label{eqpartialsumcoef}
\sum_{n = 1}^m {{C_n}(y)} &=&\sum_{n = 1}^m
{{U_n}\bigl(\sqrt y \bigr){{\bigl(\sqrt y \bigr)}^n}}\nonumber\\
&=& \sum_{n =
1}^m { \frac{\sin(n+1) \theta}{ \sin\theta} \cdot
(\cos\theta)^n } \\
& =& \frac{\cos^2 \theta}{ \sin^2 \theta} [1-(\cos\theta)^{m-1}
\cdot\cos{(m+1)\theta}]>0.\nonumber
\end{eqnarray}
Since $N$ is nonnegative and $e' \ge e'N$, we immediately get
$e'N^{n} \ge e'N^{n+1} \ge0$ for all $n \ge0$. Combining this fact
with (\ref{eqpartialsumcoef}), we obtain
%e15 #&#
%
\begin{eqnarray} \label{eqpartialvecmat}
{e'\sum_{n = 1}^k {{C_n}(y){N^n}} }
\ge\sum_{n = 1}^k {{C_n}(y)e'{N^k}} \ge0'\qquad \mbox{for all $k \ge1$}.
\end{eqnarray}
Therefore, from (\ref{eqC0}), (\ref{eqchebmatrix}) and
(\ref{eqpartialvecmat}) we conclude that, for all $y \in(0, 1),$
\begin{eqnarray*}
{e'{\bigl((1 - y)I + y{(I - N)^2}\bigr)^{ - 1}}} & = & e'\sum_{n = 0}^\infty
{{C_n}(y){N^n}} \\
& = & \lim_{k \to\infty} e'\sum_{n = 1}^k {{C_n}(y){N^n}} + e'
\\
&\ge& 0' + e' = e' >0'.
\end{eqnarray*}
This concludes the proof of (c) and we deduce that there exists a
CQLF for the pair $(-R, -R(I-pe'))$.

To prove the existence of a CQLF for the other pair $(-R,
-(I-pe')R)$, we note that $-(I-pe')R$ has the same spectrum as
$-R(I-pe')$ and the matrix product $R(I-pe')R$ has the same spectrum
as $R^2(I-pe').$ Application of Proposition~\ref{propshorten2009}
completes the proof of Theorem~\ref{thmexistcqlf}.
\end{pf}

%s5.2 #&#
\subsection{\texorpdfstring{Proof of Theorem~\protect\ref{thmnoabd}}{Proof of Theorem 2}}

In this section we prove Theorem~{\ref{thmnoabd}}. Key to the proof
is the CQLF constructed from Theorem~\ref{thmexistcqlf}.

\begin{pf} If $\alpha=0$, then from
(\ref{eqdrfitcoef}) we know that $Y$ has the piecewise linear drift
\[
b(y)= - \beta p -R\bigl(y - p (e'y)^+\bigr).
\]
By Theorem~\ref{thmexistcqlf}, there exists a CQLF
%e16 #&#
%
\begin{equation} \label{eqlyapunof}
L(y)=y'Qy,
\end{equation}
where $Q$ is a positive definite matrix such that
%e17 #&#
%e18 #&#
%
\begin{eqnarray}
Q(-R)+(-R)^{\prime}Q&<&0, \label{eqQ1R}\\
Q\bigl(-R(I-pe')\bigr)+\bigl(-(I-ep')R'\bigr)Q&\le&0. \label{eqQ1Rsig}
\end{eqnarray}

We claim that given any positive constant $C>0$, there exists a
constant $M>0$ such that if $|y|>M$,
%e19 #&#
%
\begin{equation}
\label{ineqdriftineq} (\nabla L(y) )'b(y) \le-C.
\end{equation}
We discuss the cases $e'y<0$ and $e'y \ge0$ separately.

\textit{Case} 1. $e'y < 0$. In this case, we have
\[
(\nabla L (y) )'b(y) = y'
[Q(-R)+(-R)^{\prime}Q ] y- 2 \beta p'Q y.
\]
By (\ref{eqQ1R}), the quadratic term dominates if $|y|$ is large.
Thus there exists a constant $M_1>0$ such that when $e'y < 0$ and
$|y|>M_1$,
\[
\label{ineqdrfitneg}
(\nabla L (y) )'b(y) \le
-C.
\]

\textit{Case} 2. $e'y\ge0$. In this case, we have
%e20 #&#
%
\begin{equation} \label{eqcase2}
\qquad\hspace*{8pt}(\nabla L (y) )'b(y) = y'\bigl[Q\bigl(-R(I-pe')\bigr)+\bigl(-(I-ep')R'\bigr)Q\bigr]y- 2 \beta
p'Qy.
\end{equation}
To overcome the difficulty caused by the singularity of $-R(I-pe'),$
we decompose $y$ as follows:
%e21 #&#
%
\begin{equation} \label{eqdecompy}
y=a p+\xi,
\end{equation}
where $\xi' p=0$ and $a \in\mathbb{R}$. Then we have
%e22 #&#
%
\begin{equation}
\label{eqkxi} |y|^2= |ap|^2+|\xi|^2 \quad\mbox{and}\quad
e'y=a+e'\xi\ge0.
\end{equation}
Note that $p'[Q(-R(I-pe'))+(-(I-ep')R')Q] p=0.$
Using (\ref{eqQ1Rsig}), we obtain
$p'[Q(-R(I-pe'))+(-(I-ep')R')Q]=0'$. This immediately implies
$p'[Q(-R(I-pe')]=0.$ Since $(I-pe')$ has a
simple zero eigenvalue, we have
\[
\label{eqpprimeQ}
p'Q=b e'R^{-1} \qquad\mbox{for some $b \ne0$.}
\]
Using this fact, we rewrite the left-hand side of (\ref{eqQ1Rsig})
as
%e23 #&#
%
\begin{eqnarray}\label{eqQsigrewrite}
&& Q\bigl(-R(I-pe')\bigr)+\bigl(-(I-ep')R'\bigr)Q
\nonumber
\\[-8pt]
\\[-8pt]
\nonumber
&&\qquad= \bigl((I-ep')R'\bigr)\cdot\bigl(-QR^{-1}-(R^{-1})'Q\bigr) \cdot\bigl(R(I-pe')\bigr).
\end{eqnarray}
After left-multiplying by $(R^{-1})'$ and right-multiplying by
$R^{-1}$ in (\ref{eqQ1R}), we deduce that $[-QR^{-1}-(R^{-1})'Q]$
is a negative definite matrix. Moreover, since $\xi' p=0$, from
(\ref{eqdecompy}) and (\ref{eqQsigrewrite}) we know that there
exists some $c>0$ such that
%e24 #&#
%
\begin{eqnarray} \label{eqxineg}
&&y'\bigl[Q\bigl(-R(I-pe')\bigr)+\bigl(-(I-ep')R'\bigr)Q\bigr]y
\nonumber\\
&&\qquad=
y'\bigl[(I-ep')R' \cdot
\bigl(-QR^{-1}-(R^{-1})'Q\bigr) \cdot R(I-pe')\bigr]y
\nonumber
\\[-8pt]
\\[-8pt]
\nonumber
&&\qquad= \xi'\bigl((I-ep')R'\bigr)\cdot\bigl(-QR^{-1}-(R^{-1})'Q\bigr) \cdot\bigl(R(I-pe')\bigr)\xi
\\
&&\qquad\le - c |\xi|^2.\nonumber
\end{eqnarray}
Therefore, from (\ref{eqcase2}) we have that for any $y$ with $e'y
\ge0,$
%e25 #&#
%e26 #&#
%
\begin{eqnarray}
(\nabla L(y) )'b(y)
&\le& - c |\xi|^2 - 2 \beta p'Q \xi- 2\beta a p'Qp \label
{ineqdrfitposwithk} \\
&\le& - c |\xi|^2- 2 \beta p'Q \xi+2 \beta p'Qp e'\xi,
\label{ineqdrfitpos}
\end{eqnarray}
where the second inequality is obtained from (\ref{eqkxi}),
$\beta>0$ and $p'Qp>0$. For $|y|$ large, if $|\xi| \ge r$ for some
large constant $r$, we obtain $ (\nabla L (y) )'b(y) \le-C$ since
the quadratic term $- c |\xi|^2$ in (\ref{ineqdrfitpos}) dominates.
If $|\xi| <r$ and $|y|$ large, we deduce from (\ref{eqkxi}) that
$a$ must be positive and large, that is,
\[
\label{ineqaposlarge}
a \ge\frac{1}{|p|} \sqrt{|y|^2-r^2}.
\]
Hence, the dominating term in (\ref{ineqdrfitposwithk}) is $- 2\beta
a p'Qp$ and we immediately obtain $ (\nabla L (y) )'b(y) \le-C$
whenever $|y|$ is large. Therefore, there exists a constant $M_2>0$
such that when $e'y \ge0$ and $|y|>M_2$,
\[
\label{ineqdrfitposC}
(\nabla L (y) )'b(y) \le-C.
\]
On setting $M=\max\{M_1, M_2\}$, we immediately get
(\ref{ineqdriftineq}).\vadjust{\goodbreak}

Now set $C=|\sum_{i,j} {Q_{ij} (\sigma\sigma')_{ij} }|+1.$
Equations~(\ref{eqlyapunof}) and (\ref{ineqdriftineq}) imply that
for $|y|>M$,

\[
\label{eqgene} GL(y) = \sum_{i,j} {{Q}_{ij} (\sigma\sigma')_{ij} }
+ (\nabla L (y) )'b(y) \le-1.
\]
The proof of Theorem~{\ref{thmnoabd}} is complete after applying
Proposition~\ref{propprec}.
\end{pf}

%s5.3 #&#
\subsection{\texorpdfstring{Proof of Theorem~\protect\ref{thmqedabd}}{Proof of Theorem 3}} \label{secpfofthm3}

In this section we prove Theorem~{\ref{thmqedabd}}. Throughout this
section, $C$ is a generic positive constant which may differ from
line to line but is independent of $y$.

By Theorem~\ref{thmexistcqlf}, there exists a positive definite
matrix $\widetilde Q$ with $|\widetilde Q|=1$ such that
%e27 #&#
%e28 #&#
%
\begin{eqnarray}
\widetilde Q(-R)+(-R)^{\prime} \widetilde Q&<&0, \label{eqZCQLF}\\
\widetilde Q\bigl(-(I-pe')R\bigr)+\bigl(-R'(I-ep')\bigr) \widetilde Q &\le&0.
\label{eqZCQLFsig}
\end{eqnarray}
We construct a nonquadratic Lyapunov function $V \in
C^2(\mathbb{R}^K)$ as follows. Let
%e29 #&#
%
\begin{equation}\label{eqlyapfuncV}
V(y)=(e'y)^2 + \kappa[y-p \phi(e'y)]' \widetilde Q [y-p \phi(e'y)],
\end{equation}
where $\kappa$ is a positive constant to be decided later and
$\phi(x)$ is a real-valued $C^2(\mathbb{R})$ function, approximating
$x \mapsto x^+$. Specifically, fix $\varepsilon>0$ and let
\[
\phi(x) =
\cases{ x, &\quad $\mbox{if $x \ge0$,}$
\vspace*{2pt}\cr
 - \frac{1}{2}\varepsilon, &\quad $\mbox{if $x\le-\varepsilon$,}$\vspace
*{2pt}\cr
\mbox{smooth}, & \quad$\mbox{if $ -\varepsilon< x <0 $}.$}
\]
We piece $x \ge0$ and $x\le-\varepsilon$ together in a smooth way
such that $\phi$ is in $C^2(\mathbb{R}),$ $ - \frac{1}{2}\varepsilon
\le\phi(x) \le x^+$ and $0 \le\dot\phi(x) \le1$ for any $x \in
\mathbb{R}$, where $\dot\phi$ is the derivative of $\phi$. This
function $\phi$ evidently exists. Note that $V \in
C^2(\mathbb{R}^K),$ but that it is not a CQLF due to its
nonquadratic nature. We summarize the key result in the following
proposition, which implies Theorem~{\ref{thmqedabd}}.

%pr4 #&#
%
\begin{proposition}
\label{proplyafosand} If $\alpha>0$, there exists a constant $C>0$
such that when $|y|$ is large enough, we have
%e30 #&#
%
\begin{eqnarray}
\label{ineqlya-foster}\qquad\quad  (\nabla V(y))'b(y) \le-C |y|^2 \quad
\mbox{and}\quad \biggl|\frac{{\partial^2 V}}{{\partial y_i\,
\partial y_j }}(y)\biggr| \le C |y| \qquad\mbox{for any $i, j$.}
\end{eqnarray}
Consequently, when $|y|$ is large,
%e31 #&#
%
\begin{equation} \label{eqGVy}
GV(y) \le-C |y|^2 \le-1.
\end{equation}
\end{proposition}

\begin{pf}
We first study $(\nabla V(y))'b(y).$ From (\ref{eqlyapfuncV}), we
have for all \mbox{$y \in\mathbb{R}^K,$}
%e32 #&#
%
\begin{equation}
\label{eqVdir} (\nabla V(y))'=2(e'y)e' + 2\kappa\bigl(y'-p' \phi(e'y)\bigr)
\widetilde Q [I-pe' \dot\phi(e'y)].
\end{equation}
We discuss the cases $e'y \ge0$, $e'y \le-\varepsilon$ and
$-\varepsilon< e'y <0$ separately.\vadjust{\goodbreak}

\textit{Case} 1. $e'y \ge0$. In this case, let $x=e'y$ and
$z=y-px=(I-pe')y$, then we have
\begin{eqnarray*}
&&(\nabla V(y))'b(y)\\
&&\qquad= [2(e'y) e'+2 \kappa y'(I-ep')\tilde
Q(I-pe')]\bigl(-R(I-pe')y-\alpha
p e'y - \beta p\bigr)\\
&&\qquad= - 2\alpha{x^2} - \kappa z'[\widetilde Q(I-pe')R+
R'(I-ep')\widetilde Q]z - 2x\beta- 2xe'Rz.
\end{eqnarray*}

Suppose we have shown that there exists $C>0$ such that
%e33 #&#
%
\begin{equation} \label{eqzneg}
z'[\widetilde Q(I-pe')R+ R'(I-ep')\widetilde Q]z \ge C |z|^2,
\end{equation}
we then obtain that
\[
(\nabla V(y))'b(y) \le- 2\alpha{x^2} - \kappa C |z|^2 - 2x\beta-
2xe'Rz.
\]
Since $\alpha>0,$ we can select $\kappa>0$ large so that
$\frac{1}{2}(2\alpha{x^2} + \kappa C |z|^2)> 2|xe'Rz|$ for any $(x,
z),$ where $\kappa$ is independent of $(x,z)$ or $y$. Then we have,
\[
(\nabla V(y))'b(y) \le- \alpha{x^2} - \tfrac{1}{2}\kappa C |z|^2 -
2x\beta.
\]
Note that $|y|=|px+z| \le C|(x,z)|,$ so that $|(x,z)|$ is large
whenever $|y|$ is large. We conclude that for $|y|$ large,
\begin{eqnarray*}
(\nabla V(y))'b(y) &\le& -C|(x,z)|^2 \\
&\le& -C|y|^2.
\end{eqnarray*}

It remains to prove (\ref{eqzneg}). We use a similar argument as
for (\ref{eqxineg}). Observe that
\[
\label{eqtildeQsigdir}
(R^{-1}p)'[\widetilde Q(I-pe')R+ R'(I-ep')\widetilde Q](R^{-1}p)=0,
\]
which implies that $\widetilde QR^{-1}p=b e$ for some $b \in
\mathbb{R}.$ Thus, we obtain
%e34 #&#
%
\begin{eqnarray}\label{eqzpos}
&&z'[\widetilde Q(I-pe')R+ R'(I-ep')\widetilde Q]z
\nonumber
\\[-8pt]
\\[-8pt]
\nonumber
&&\qquad=z'R'(I-ep')[(R^{-1})'\widetilde Q+ \widetilde QR^{-1}] (I-pe')Rz.
\end{eqnarray}
Since $R$ is a nonsingular M-matrix, $R^{-1}$ is a nonnegative
matrix Berman and Plemmons [(\citeyear{bermanplemmons}), Chapter~6], and we deduce that
%e35 #&#
%
\begin{equation} \label{eRinvp}
e'R^{-1}p>0.
\end{equation}
This implies that $(I-pe')Rz \ne0$ since $e'z=e'(I-pe')y=0$ in this
case. From~(\ref{eqZCQLF}) we know that $(R^{-1})'\widetilde Q+
\widetilde QR^{-1}$ is a positive definite matrix. Now~(\ref{eqzneg}) follows from (\ref{eqzpos}).

\textit{Case} 2. $e'y < -\varepsilon$. In this case, we have
$\phi(e'y)=-\frac{1}{2}\varepsilon$ and $\dot\phi(e'y)=0$. From~(\ref{eqZCQLF}), there exists $C>0$ such that
\begin{eqnarray*}
(\nabla V(y))'b(y) &=&\bigl(2(e'y) e'+ 2\kappa y' \widetilde Q+ \kappa
\varepsilon p' \widetilde Q\bigr)(-Ry- \beta
p)\\
&=&- 2\kappa\bigl[y'(\widetilde QR+R'\widetilde Q)y+ \tfrac{1}{2}
(\varepsilon p' \widetilde Q R + \beta p'\widetilde Q)
y+ \tfrac{1}{2} \varepsilon\beta p'\widetilde Qp \bigr] \\
&&{}-2e'y \cdot(e'Ry+ \beta)\\
&\le& - 2\kappa\bigl[C |y|^2+ \tfrac{1}{2}(\varepsilon p' \widetilde
Q R + \beta
p'\widetilde Q)y+ \tfrac{1}{2} \varepsilon\beta p'\widetilde Qp \bigr]
\\
&&{}-2e'y \cdot(e'Ry+ \beta)\\
&\le& - 2\kappa\bigl[C |y|^2+ \tfrac{1}{2} (\varepsilon p' \widetilde
Q R + \beta
p'\widetilde Q)y+ \tfrac{1}{2} \varepsilon\beta p'\widetilde Qp \bigr]\\
&&{}+
\kappa C (|y|^2 + |y|) \\
&\le& -\kappa(C |y|^2-C |y|-C),
\end{eqnarray*}
where $\kappa$ is again chosen to be independent of $y$, but large
enough such that $|2e'y \cdot(e'Ry+ \beta)|< \kappa C (|y|^2+|y|).$
Thus for $|y|$ large and $e'y <-\varepsilon$, we have
\[
(\nabla V(y))'b(y) \le -C |y|^2.
\]

\textit{Case} 3. $-\varepsilon\le e'y \le0 $. In this case we use the
property that $0\le\dot\phi(e'y)\le1$. Note that we have
\begin{eqnarray*}
&&(\nabla V(y))'b(y) \\
&&\qquad=\bigl(2(e'y)e' + 2 \kappa\bigl(y'-p' \phi(e'y)\bigr) \widetilde Q \bigl(I-pe' \dot
\phi(e'y)\bigr) \bigr) (-R y-
\beta p)\\
&&\qquad= 2 e'y e'(-Ry-\beta p) \\
&&\qquad\quad{}+ 2 \kappa\dot\phi(e'y)\bigl (y'-p' \phi(e'y)\bigr) \widetilde Q(I-pe')
(-Ry- \beta p) \\
&&\quad\qquad{}+ 2 \kappa\bigl(1-\dot\phi(e'y)\bigr)\bigl(y'-p' \phi(e'y)\bigr) \widetilde Q (-R
y- \beta p).
\end{eqnarray*}
We write
\[
y= a R^{-1} p + \xi,
\]
where $\xi$ is orthogonal to
$R^{-1} p$ and $a \in\mathbb{R}$, so that
%e36 #&#
%
\begin{equation} \label{eqydecomdis}
|y|^2=c a^2 + |\xi|^2 \qquad\mbox{for some $c>0.$}
\end{equation}
From (\ref{eRinvp}), we have $e'R^{-1}p>0.$ Without loss of
generality we assume that $e'R^{-1}p=1$. Then $e'y= a +e'\xi$ and
we get
%e37 #&#
%
\begin{eqnarray}
\label{ineqcase3}
&&(\nabla V(y))'b(y) \nonumber\\
&&\qquad= -2 (a+e'\xi) ( \beta+e'R \xi+ a)
\nonumber
\\
&&\quad\qquad{} + \kappa\dot\phi(e'y) \bigl(\xi'\bigl[\widetilde
Q\bigl(-(I-pe')R\bigr)+\bigl(-(I-pe')R\bigr)' \tilde Q\bigr]
\xi
\nonumber
\\[-8pt]
\\[-8pt]
\nonumber
&&\hspace*{118pt}\qquad\quad{}-2 p' \widetilde Q(I-pe')R\xi\phi(e'y)\bigr)\\
&&\qquad\quad{} +\kappa\bigl(1-\dot\phi(e'y)\bigr)\nonumber\\
&&\qquad\phantom{+}\quad{}\times\bigl (y'[-\widetilde Q R-R' \widetilde Q]y
+ \beta y' \widetilde Q p- \phi(e'y) p' \widetilde Q R y
-p'\widetilde Q p \beta\bigr).\nonumber
\end{eqnarray}
Since $\xi'R^{-1}p=0$, one checks, as for (\ref{eqzneg}), that there
exists a constant $C>0$ such that
%e38 #&#
%
\begin{equation} \label{eqxineg2}
\xi'\bigl[\widetilde Q\bigl(-(I-pe')R\bigr)+\bigl(-(I-pe')R\bigr)'\widetilde Q \bigr] \xi\le-C
|\xi|^2.
\end{equation}
Moreover, from (\ref{eqZCQLF}) and (\ref{eqydecomdis}), we deduce
that
%e39 #&#
%
\begin{equation} \label{eqzcqlfneg2}
y'[-\widetilde Q R-R' \widetilde Q]y \le-C |y|^2= -C a^2 - C
|\xi|^2.
\end{equation}
Substituting (\ref{eqxineg2}) and (\ref{eqzcqlfneg2}) into
(\ref{ineqcase3}), and using $0\le\dot\phi(e'y)\le1$ as well as
$|\phi(e'y)| \le\varepsilon$, we obtain
%e40 #&#
%
\begin{eqnarray}\label{eqkappaCxi}
&&{(\nabla V(y))'b(y) }
\nonumber
\\[-8pt]
\\[-8pt]
\nonumber
&&\qquad\le-2 (a^2 + C |a| |\xi| + C|a|) + \kappa(-C
|\xi|^2+ C|\xi|+C|a|+C).
\end{eqnarray}
Since $e'y= a +e'\xi\in[-\varepsilon, 0],$ we must have $|a| \le C+
|\xi|$ and consequently $|y| \le C |a|+ |\xi| \le C |\xi|+C$. Thus
for $|y|$ large, we can choose $\kappa$ large so that the dominating
term in (\ref{eqkappaCxi}) is $-\kappa C |\xi|^2$. Using the fact
that $|y|^2 \le C |\xi|^2$ when $|y|$ is large, we then deduce that
there exists a constant $C>0$ such that for $|y|$ large,
\[
(\nabla V(y))'b(y) \le-C |y|^2.
\]
This concludes the proof for the third case.

On combining the above three cases we obtain that, for $|y|$ large,
\[
(\nabla V(y))'b(y) \le-C |y|^2,
\]
as claimed in the proposition.

We now proceed to study the second derivative of $V$, which is
denoted by~$\ddot V.$ We also write $\ddot\phi$ for the second
derivative of $\phi.$ From (\ref{eqVdir}), we find
%e41 #&#
%
\begin{eqnarray}
\label{eqVseconddir}
\ddot V(y)&=& 2ee'+ 2 \kappa\bigl[\widetilde Q+
ee'\cdot p'\widetilde Q p \bigl(\ddot\phi(e'y)\phi(e'y)+\dot\phi
(e'y)^2\bigr)
\nonumber
\\[-8pt]
\\[-8pt]
\nonumber
&&\hspace*{46pt}{}- (\widetilde Qpe'+ ep'\widetilde Q) \dot\phi(e'y)- ee' \cdot
y'\widetilde Qp \ddot\phi(e'y) \bigr].
\end{eqnarray}
If $e'y \notin[-\varepsilon, 0]$, we obtain $0\le\dot\phi(e'y)\le
1$ and $\ddot\phi(e'y)=0$. Therefore, for any $i, j$, there exists
some $C>0$ such that
\[
\biggl|\frac{{\partial^2 V}}{{\partial y_i\,
\partial y_j }}(y)\biggr| \le C.
\]
If $e'y \in[-\varepsilon, 0]$, then $|\ddot\phi(e'y)| \le C$ for
some $C>0$ since $\phi\in C^2(\mathbb{R})$ and $[-\varepsilon, 0]$ is
compact. Moreover, since $0\le\dot\phi(e'y)\le1,$ the dominating
term in (\ref{eqVseconddir}) is $ -2 \kappa ee' \cdot y'\widetilde
Qp \ddot\phi(e'y)$ for $|y|$ large. This implies that if $e'y \in
[-\varepsilon, 0]$ and $|y|$ is large, then there exists a constant
$C>0$ such that for any $i, j$,
\[
\biggl|\frac{{\partial^2 V}}{{\partial y_i\,
\partial y_j }}(y)\biggr| \le C|y|,
\]
where $C$ is independent of $y$. This concludes the proof of
(\ref{ineqlya-foster}). Now for $|y|$ large, we deduce from
(\ref{ineqlya-foster}) that
\[
GV(y)= (\nabla V(y))'b(y) + \frac{1}{2}\sum_{i,j} {(\sigma
\sigma')_{ij} \frac{{\partial^2 V}}{{\partial y_i\,
\partial y_j }}(y)} \le-C |y|^2 \le-1.
\]
The proof of Proposition~{\ref{proplyafosand}} is complete.
\end{pf}

\begin{pf*}{Proof of Theorem~\ref{thmqedabd}}
In order to show that $Y$ is positive recurrent and has a unique
stationary distribution, we only have to check that $V(y) \to
\infty$ as $ |y| \to\infty$ in view of Proposition~\ref{propprec}
and (\ref{eqGVy}).

Let $x=e'y$ and $z=y-px^+$, then $|y|^2 \le C (x^2 +|z|^2)$. We can
rewrite (\ref{eqlyapfuncV}) as follows:
\begin{eqnarray*}
V(y)&=&x^2+ \kappa\bigl(y'-p' \phi(x)\bigr) \widetilde Q \bigl(y-p \phi(x)\bigr)
\nonumber\\
&\ge& x^2+ C |y-p \phi(x)|^2 \nonumber\\
&=& x^2+ C \bigl|z+ p\bigl(x^+ - \phi(x)\bigr)\bigr|^2 \nonumber\\
&\ge& x^2 + C|z|^2 - C \varepsilon^2 \nonumber\\
&\ge& C|y|^2 - C \varepsilon^2,
\end{eqnarray*}
where the second last inequality uses the fact $ 0\le x^+ -
\phi(x)\le\frac{1}{2}\varepsilon. $ Therefore, $V(y) \to\infty$
as $
|y| \to\infty$ and we conclude that $Y$ has a unique stationary
distribution.

To prove that $Y$ is exponentially ergodic, we observe from
(\ref{eqlyapfuncV}) that there exists some $C>0$ such that $V(y)\le
C |y|^2+C$ for all $y \in\mathbb{R}^K$. Moreover, (\ref{eqGVy})
implies that for $|y|$ large,
\[
GV(y) \le-C V(y)+C.
\]
Putting this together with the fact that $V\in C^2 (\mathbb{R}^K)$,
we know that there exist some $c>0$ and $d< \infty$ such that
\[
GV(y) \le-c V(y)+d\qquad \mbox{for any $y \in\mathbb{R}^K$}.
\]
Since $V \ge0,$ Proposition~\ref{propexpergo} implies that $Y$ is
$f$-exponentially ergodic, where $f=V+1$. In particular, $Y$ is
exponentially ergodic since $f \ge1.$
\end{pf*}

\begin{appendix}
%s6 #&#
\section{\texorpdfstring{Proof of Proposition~\lowercase{\protect\ref{propshorten2009}}}{Proof of Proposition 3}}\label
{appendixprop3pf}
We first outline the key idea behind the proof.
Suppose that $(B, g)$ is not controllable or that $(B, h)$ is not
observable in the CQLF existence problem. Then we can ``reduce'' them
to suitable subspaces such that $(B_1, g_1)$ is controllable and
$(B_1, h_1)$ is observable, where $B_1$ is a new matrix of lower
dimension than $B$ and similarly for $g_1, h_1.$ In the process of
``reduction,'' two desired properties are preserved: (a) $B(B-gh')$
has no real negative eigenvalues if and only if $B_1(B_1-g_1 h_1')$
has no real negative eigenvalues; (b) $(B, B-gh')$ has a CQLF if and
only if $(B_1, B_1-g_1h_1' )$ has a CQLF. Therefore, applying
Theorem~3.1 in \citet{Shorten09} to $(B_1, B_1 - g_1 h_1')$ yields
the result.

To make the ideas concrete, we now introduce a lemma
giving an equivalent formulation of the CQLF existence problem,
which makes the ``reduction'' possible. The lemma is an analog of
Proposition 2 in \citet{King06}. In \citet{King06}, each matrix of the
pair is nonsingular while in our case one of the matrices is
singular.

%le1 #&#
%
\begin{lemma}
\label{lemdualcqlf} Suppose that all eigenvalues of the matrix $B$
have negative real part and all eigenvalues of $B-gh'$ have negative
real part, except for a simple zero eigenvalue. Then the following
statements are equivalent:
\begin{longlist}[(a)]

\item[(a)] The pair $(B, B-gh')$ does not have a CQLF.

\item[(b)] There are positive semidefinite matrices $X$ and $Z$ such
that
\begin{eqnarray*}
&&BX+XB'+(B-gh')Z+Z(B'-hg')=0, \\
&& BX+XB' \ne0\quad \mbox{and}\quad (B-gh')Z+Z(B'-hg')\ne0.
\end{eqnarray*}

\item[(c)] There are nonzero, positive semidefinite matrices $X$ and
$Z$ such that
%e42 #&#
%
\begin{equation}
\label{eqdualcqlf}
BX+XB'+(B-gh')Z+Z(B'-hg')=0,
\end{equation}
where $Z \ne c B^{-1}g g'(B^{-1})'$ for any $c \in\mathbb{R}.$
\end{longlist}
\end{lemma}

\begin{pf}
We first prove the equivalence of (a) and (b). To set up the
notation, let $S^{K \times K}$ be the space of real symmetric $K
\times K$ matrices. For an arbitrary matrix $A \in\mathbb{R}^{K
\times K}$, define the linear operator $L_A $ on $S^{K \times K}$ by
%e43 #&#
%
\begin{equation} \label{eqlyaoperator}
L_A\dvtx S^{K \times K} \rightarrow S^{K \times K},\qquad
{L_A(H)=AH+HA'}.
\end{equation}
It is well known that if $A$ has eigenvalues $\{\lambda_i\} $ with
eigenvectors $\{v _i\} $, then $L_{A}$ has eigenvalues $\{
\lambda_i + \lambda_j \}$ with eigenvectors $\{v_i v_j' + v_j
v_i'\}$ for all $i \le j$. Since all eigenvalues of the matrix $B$
have negative real part, $L_B$ is invertible.

Following \citet{King06}, we formulate the CQLF existence problem in
terms of separating convex cones in $S^{K \times K}$. Define $
\Cone(B) = \{ L_B (X)|X \ge0\} $ and $ \Cone(B - gh') = \{ L_{(B -
gh')}(Z) |Z \ge0\} $. Both are closed convex cones in $S^{K \times
K}$. Let $S^{K \times K}$ be equipped with the usual Hilbert--Schmidt
inner product $\langle X, Z \rangle= \tr(XZ)$. We obtain that for
any $Q \in S^{K \times K}$,
\[
\langle X, QB+B'Q \rangle=\langle Q, BX+XB' \rangle= \langle Q,
L_B(X) \rangle.
\]
Note that for a nonzero positive semidefinite matrix $X$,
we have $QB+B'Q<0$ if and only if $\langle X, QB+B'Q
\rangle<0,$ where the ``if'' part can be checked by taking $X=xx'$
for any nonzero $x \in\R^K$, and the ``only if'' part follows from
the spectral decomposition of the positive semidefinite matrix $X.$
Therefore, we have $QB+B'Q<0$ if and only if $\langle Q, M \rangle
<0 $ for all nonzero $M \in\Cone(B)$. Using a similar argument one
finds that $Q(B-gh')+(B-hg')Q \le0$ if and only if $\langle
Q,T\rangle\le0 $ for all nonzero $T \in\Cone(B-gh')$. Moreover,
since $B$ only has eigenvalues with negative real part, we deduce
that $QB+B'Q<0$ for $Q \in S^{K \times K}$ implies that $Q$ is
positive definite by Theorem 2.2.3 in \citet{hornandjohnson94}. By
definition of CQLF, we thus obtain that $(B, B-gh')$ has a CQLF if
and only if there exists a $Q \in S^{K \times K} $ such that
$QB+B'Q<0$ and $Q(B-gh')+(B-hg')Q \le0.$ Equivalently, $(B, B-gh')$
has a CQLF if and only if there exists a $Q \in S^{K \times K}$ such
that $\langle Q,M \rangle>0 $ for all nonzero $M \in\Cone(-B)$ and
$\langle Q,T\rangle\le0 $ for all nonzero $T \in\Cone(B-gh').$
Therefore, finding a CQLF for the pair $(B, B-gh')$ is the same as
finding a separating hyperplane in $S^{K \times K}$ for $\Cone( -
B)$ and $\Cone(B - gh').$ By the separating hyperplane theorem, we
conclude that $(B, B-gh')$ not having a CQLF is equivalent to
$\Cone( - B)$ and $\Cone(B - gh')$ having nonzero intersection. This
completes the proof of the equivalence of (a) and (b).

We now turn to the equivalence of (b) and (c), for which we use
the aforementioned spectral properties of the linear operator
(\ref{eqlyaoperator}). Since $L_B$ is invertible, we deduce that
$L_B(X)=0$ is equivalent to $X=0$. We know that all eigenvalues of
$(B-gh')$ have negative real part except for a simple zero
eigenvalue, hence, $L_{(B - gh')}$ also has a simple zero eigenvalue
with eigenvector $cB^{-1}gg' (B^{-1})'$ for some nonzero $c \in
\mathbb{R}$ while all of its other eigenvalues have negative real
part. Consequently, $(B-gh')Z+Z(B-gh')' \ne0$ is equivalent to $Z
\ne c B^{-1}g g'(B^{-1})'$ for any $c \in\mathbb{R}.$ The proof of
the lemma is complete.~%
\end{pf}

\begin{pf*}{Proof of Proposition~\ref{propshorten2009}} In view of Theorem 3.1
of \citet{Shorten09}, we need to check that controllability of
$(B,g)$ and observability of $(B,h)$ need not be verified in the
CQLF existence problem. Recall that controllability of $(B,g)$ means
that the vectors $g,Bg,B^2 g, \ldots$ span $\mathbb{R}^K$, and
observability of $(B,h)$ means that the vectors $h,B'h,(B')^2 h,
\ldots$ span $\mathbb{R}^K.$ To simplify the notation, let
$\widetilde B =B-gh'.$

We first show that in the CQLF existence problem for the pair $(B,
B-gh')$, we can assume without loss of generality that $(B,g)$ is
controllable. The proof relies on Lemma~\ref{lemdualcqlf}. Let $U$
be the span of vectors $g,Bg,B^2 g \ldots.$ Suppose $U$ is a proper
subspace of $\mathbb{R}^K$ with $\dim(U)<K$, and note that
$\mathbb{R}^K = U \oplus U^ \bot$ where $U^ \bot$ is the orthogonal
complement of $U$. In view of this decomposition, we perform a
change of basis and rewrite $B$, $\tilde B$, $X$ and $Y$ in the
block form
%e44 #&#
%
\begin{eqnarray}\label{eqrewriteblock}
B &=& \pmatrix{
{{B_1}} & {{B_2}} \vspace*{2pt}\cr
0 & {{B_3}} }
,\qquad \widetilde B = \pmatrix{
{{{\widetilde B}_1}} & {{{\widetilde B}_2}} \vspace*{2pt}\cr
0 & {{B_3}} },
\nonumber
\\[-8pt]
\\[-8pt]
\nonumber
X &=& \pmatrix{
{{X_1}} & {{X_2}} \vspace*{2pt}\cr
{{X_2'}} & {{X_3}} }
,\qquad
 Z = \pmatrix{
{{Z_1}} & {{Z_2}} \vspace*{2pt}\cr
{{Z_2'}} & {{Z_3}} },
\end{eqnarray}
where $B-\widetilde B= gh'$ and $g, h$ are represented in the new
basis. We use the same notation for the matrices and vectors after
the change of basis to save space, and we remark that the orthogonal
transformation does not affect the existence of a CQLF for the pair
$(B, \tilde B)$ or the existence of real negative eigenvalues of $B
\tilde B$. Namely, for any orthonormal matrix $O \in\mathbb{R}^K$,
one readily checks that the pair $(B, \tilde B)$ has a CQLF if and
only if the pair $(OBO', O \tilde B O')$ has a CQLF. Furthermore, $B
\tilde B$ has no real negative eigenvalues if and only if $(OBO')(O
\tilde B O')$ has no real negative eigenvalues. Let $g_1, h_1$ be
the orthogonal projection of $g, h$ on the subspace $U,$ so that
$B_1-\widetilde B_1= g_1 h_1'$. Since $U$ is the span of the vectors
$g,Bg,B^2 g \ldots,$ we deduce that $g_1,B_1g_1,B_1^2 g_1 \ldots$
span $U$ by (\ref{eqrewriteblock}), that is, $(B_1, g_1)$ is
controllable. We now use Lemma~{\ref{lemdualcqlf}} to argue that
there exists a CQLF for $(B, \widetilde B)$ if and only if there
exists a CQLF for $(B_1,\widetilde B_1 )$, where $(B_1, g_1)$ is
controllable. Note that (\ref{eqrewriteblock}) implies, using
(\ref{eqdualcqlf}) in Lemma~\ref{lemdualcqlf},
\[
B_3(X_3+Z_3)+(X_3+Z_3)B_3'=0.
\]
Equivalently,
\[
L_{B_3}(X_3+Z_3)=0,
\]
where the linear operator $L_{B_3}$ is defined in
(\ref{eqlyaoperator}). Since $B$ has only eigenvalues with
negative real part, $B_3$ also has this property.
This implies the linear operator $L_{B_3}$ is
invertible. We thus obtain $X_3+Z_3=0.$ Using the fact that $X$
and $Z$ are positive semidefinite, we deduce that $X_3=Z_3=0,$ and
consequently \mbox{$X_2=Z_2=0$}. This leads to
%e45 #&#
%
\begin{equation}
\label{eqsub} B_1X_1+X_1{B_1'}+{\widetilde B_1}
Z_1+Z_1{\widetilde B_1'}=0.
\end{equation}
Thus, for the pair $(B, B-gh')$, the existence of nonzero $X, Z \ge
0$ such that~(\ref{eqdualcqlf}) holds implies the existence of
nonzero $X_1, Z_1 \ge0 $ such that~(\ref{eqsub}) holds. Conversely,
if there exists nonzero $X_1, Z_1 \ge0 $ such that (\ref{eqsub})
holds, setting $X_2=X_3=Z_2=Z_3=0,$ we then obtain that there exists
nonzero $X, Z \ge0$ such that~(\ref{eqdualcqlf}) holds. Since
$B-gh'$ has only eigenvalues with negative real part except for a
simple zero eigenvalue, so does $B_1-g_1h_1'.$ For $c \in
\mathbb{R}$, since $g\in U,$ one finds that $g'(B^{-1})' =
(g_1'(B_1^{-1})', 0')$ by (\ref{eqrewriteblock}). Thus $Z \ne c
B^{-1}g g'(B^{-1})'$ is equivalent\vadjust{\goodbreak} to $Z_1 \ne c B_1^{-1}g_1
g_1'(B_1^{-1})'.$ Putting these\vspace*{2pt} facts together, we apply
Lemma~\ref{lemdualcqlf} to conclude that $(B, \widetilde B)$ has no
CQLF if and only if $(B_1,\widetilde B_1 )$ has no CQLF, where
$(B_1, g_1)$ is controllable. Therefore, without loss of generality,
we can assume that $(B,g)$ is controllable in the CQLF existence
problem for the pair $(B, B-gh')$.

We next show that without loss of generality we can assume that $(B,
h)$ is observable in the CQLF existence problem for the pair $(B,
B-gh').$ Note that for $Q>0,$ we have $QB+B'Q<0$ and
$Q(B-gh')+(B'-hg')Q \le0$ if and only if $Q^{-1}B'+B Q^{-1}<0$ and
$Q^{-1}(B-hg')+(B'-gh')Q^{-1} \le0.$ Hence, $(B, B-gh')$ has a CQLF
if and only if $(B', B'-hg')$ has a CQLF. From the preceding
paragraph, we know that in the CQLF existence problem for the pair
$(B', B'-hg')$, we can assume that $(B', h)$ is controllable without
loss of generality. By definition, $(B', h)$ being controllable is
the same as $(B, h)$ being observable. Therefore, we conclude that
we can assume without loss of generality that $(B, h)$ is
observable.

Finally, we argue that the pair $(B, B-gh')$ has a CQLF if and only
if the matrix product $B(B-gh')$ has no real negative eigenvalues.
Assuming that $(B, g)$ is controllable and that $(B, h)$ is
observable, Theorem~3.1 in \citet{Shorten09} states that $(B, B-gh')$
has a CQLF if and only if the matrix product $B(B-gh')$ has no real
negative eigenvalues. We have shown that we can always assume that
$(B, g)$ is controllable and that $(B, h)$ is observable in the CQLF
existence problem by reduction to proper subspaces. So it only
remains to check that in the process of reduction, the spectral
property of having no real negative eigenvalues of the matrix
product is preserved. Specifically, in the above proof that
controllability of $(B,g)$ can be assumed without loss of
generality, we obtain that $(B, B-gh')$ has a CQLF if and only if
$(B_1, B_1-g_1h_1' )$ has a CQLF, where $(B_1, g_1)$ is
controllable. We next prove that $B (B-gh')$ has no real negative
eigenvalues if and only if $B_1 (B_1-g_1 h_1')$ has no real negative
eigenvalues, that is, the desired spectral property of the matrix
product is preserved in the process of reduction from $(B, B-gh')$
to $(B_1, B_1-g_1 h_1' ).$ Observe that the spectrum of $B (B-gh') $
is the union of the spectrum of $B_1 (B_1-g_1h_1')$ and ${B_3^2}$ by
(\ref{eqrewriteblock}). Since all eigenvalues of $B_3$ have
negative real part, we deduce that $B_1 (B_1-g_1 h_1')$ having no
real negative eigenvalues is equivalent to $B (B-gh')$ having no
real negative eigenvalues. A similar argument applies for
observability instead of controllability. We have therefore
completed the proof of Proposition~\ref{propshorten2009}.\vspace*{-3pt}
\end{pf*}

%s7 #&#
\section{\texorpdfstring{Any quadratic function fails for $\alpha>0$}{Any quadratic function fails for alpha > 0}}
\label{appendixcountereg} In this section, we give a simple example
showing that, in general, no quadratic function can serve as a
Lyapunov function in the Foster--Lyapunov criterion to prove positive
recurrence of the piecewise OU process $Y$ for $\alpha>0.$ We first
introduce a lemma which implies that the matrix $-R(I-pe')-\alpha
pe'$ is nonsingular for $\alpha>0.$\vadjust{\goodbreak}
%le2 #&#
%
\begin{lemma}\label{lemalphanonsig}
If $\alpha> 0$, then all eigenvalues of the matrix
$-R(I-pe^{\prime})-\alpha pe^{\prime}$ have negative real part.
\end{lemma}

\begin{pf}
It is clear that the matrix has an eigenvalue $-\alpha$ with right
eigenvector~$p$. Suppose $\lambda\ne-\alpha$ is an eigenvalue of
the matrix with left eigenvector~$\theta$, that is,
%e46 #&#
%
\begin{equation}\label{eqlem2eq}
\theta'\bigl( -R(I - pe') - \alpha pe'\bigr) = \lambda\theta',
\end{equation}
then we obtain that $\theta' p = 0.$ It follows from
(\ref{eqlem2eq}) that $\lambda$ is an eigenvalue of the matrix $
-R(I - pe').$ Moreover, $\lambda$ cannot be zero since otherwise
$\theta'=c e'R^{-1}$ for some nonzero $c \in\mathbb{R},$
which follows from the fact that $R(I-pe')$ has a
simple zero eigenvalue. This contradicts the fact that $e'R^{-1}p
>0$ as seen in (\ref{eRinvp}). From condition (b) in the proof of
Theorem~\ref{thmexistcqlf}, we know that all nonzero eigenvalues of
the matrix $ -R(I - pe')$ have negative real part. This completes
the proof of the lemma.
\end{pf}

%le3 #&#
%
\begin{lemma} \label{lemweakbad}
Suppose that $Q$ is a real $K \times K$ positive semidefinite matrix
such that at least one of the matrices $Q(-R)+(-R')Q $ and
$Q(-R(I-pe')-\alpha pe')+ (-(I-ep')R'-\alpha ep')Q $ fails to be
negative definite. Let the quadratic function $L$ be given by
$L(y)=y'Qy$ for $y \in\mathbb{R}^K.$ Then there exists some $\beta
\in\mathbb{R}$ and $v \in\mathbb{R}^K$ such that $GL(tv) \ge0$
for any $t \ge0$. %\in\mathbb{R}$.
\end{lemma}

\begin{pf}
Suppose that $Q(-R)+(-R')Q$ fails to be negative definite, then
there exists some $\lambda\ge0$ and nonzero vector $v \in
\mathbb{R}^K$ such that $[Q(-R)+(-R')Q] v = \lambda v$ and $e'v \le
0.$ By definition of generator of $Y$ in (\ref{eqgenerator}), we
thus obtain
%e47 #&#
%
\begin{eqnarray} \label{eqGLtv}
GL(tv)& = &\sum_{i,j} {{Q}_{ij} (\sigma\sigma') _{ij} } +
(\nabla L
(tv) )'b(tv) \nonumber\\
& = & \sum_{i,j} {{Q}_{ij} (\sigma\sigma') _{ij} } + t^2 v'
[Q(-R)+(-R')Q]v -2t \beta p'Qv \\
& = & \sum_{i,j} {{Q}_{ij} (\sigma\sigma') _{ij} }+ \lambda v'v
t^2 -2t \beta p'Qv.\nonumber
\end{eqnarray}
Since $Q$ is positive semidefinite, we infer that $\sum_{i,j}
{{Q}_{ij} (\sigma\sigma') _{ij} } = \tr(Q \sigma\sigma')=
\tr(\sigma'Q \sigma) \ge0.$ Set $\beta=0.$ We conclude from
(\ref{eqGLtv}) that $GL(tv) \ge0$ for any $t \ge0$. %\in\mathbb{R}.$
A similar argument applies to the case where $Q(-R(I-pe')-\alpha pe')+
(-(I-ep')R'-\alpha ep')Q $ fails to be negative definite. The proof
of the lemma is complete.
\end{pf}

In view of Lemmas~\ref{lemalphanonsig} and~\ref{lemweakbad}, we
give the following definition of strong CQLF which is slightly
different than Definition~\ref{defcqlf} given in
Section~\ref{seccqlfdef}. For more details, refer to
\citet{ShortenNarendra03} and \citet{King06}.

%de6 #&#
%
\begin{definition}[(Strong CQLF)]
Let $A$ and $B$ be real $K \times K$ matrices having only
eigenvalues with negative real part. For $Q\in\mathbb{R}^{K\times
K},$ the quadratic form $L$ given by $L(y)=y'Qy$ for $y \in
\mathbb{R}^K$ is called a strong common quadratic Lyapunov function
$($strong CQLF$)$ for the pair $(A, B)$ if $Q$ is positive definite and
\begin{eqnarray*}
QA+A'Q&<&0, \\
QB+ B'Q&<& 0.
\end{eqnarray*}
\end{definition}

We remark that it suffices to require $Q$ to be a symmetric matrix
in the above definition by Theorem~2.2.3 in \citet{hornandjohnson94}.

We now formulate an example showing that, in general, no quadratic
function can serve as a Lyapunov function in the Foster--Lyapunov
criterion to prove positive recurrence of the piecewise OU process
$Y$ for $\alpha>0.$ Let $R$ be a matrix given by
\[
R = \pmatrix{
1 & { - 1} & 0 \vspace*{2pt}\cr
0 & 1 & { - 1} \vspace*{2pt}\cr
0 & 0 & 1 }
,
\]
so that $R$ is a nonsingular M-matrix. Let $\alpha=133$ and
$p'=[0,0,1].$

%le4 #&#
%
\begin{lemma}
For any quadratic function $L$ given by $L(y)=y'Qy$ for some real $K
\times K$ positive semidefinite matrix $Q$ and all $y \in
\mathbb{R}^K,$ there exists some $\beta\in\mathbb{R}$ and $v \in
\mathbb{R}^K$ such that $GL(tv) \ge0$ for any $t \in\mathbb{R}$ in
the above example.
\end{lemma}

\begin{pf}
In view of Lemma~\ref{lemweakbad}, it suffices to prove that there
is no strong CQLF for the pair $(-R, -R(I-pe')-\alpha pe')$ for
$\alpha>0.$ Equivalently, it suffices to show that the matrix
product $R (R(I-pe')+\alpha pe')$ has real negative eigenvalues by
Theorem~1 in \citet{King06}. One readily checks that $R
(R(I-pe')+\alpha pe')$ has three different eigenvalues: $-7$,
$5-\sqrt{82}$ and $5+\sqrt{82}$. Thus, it has two real negative
eigenvalues and we deduce that $(-R, -R(I-pe')-\alpha pe')$ has no
strong CQLF in this example. Application of Lemma~\ref{lemweakbad}
completes the proof of the lemma.
\end{pf}
\end{appendix}

\section*{Acknowledgments}
We thank Jim Dai for encouraging us
throughout the project and for detailed comments. We are also
grateful to Amarjit Budhiraja for helpful discussions on exponential
ergodicity, and to two anonymous referees for their thoughtful
comments.

% imsref loaded by akundreckaite, 2012-06-12 08:58:30
%

%

%suskaldyti doi

\printaddresses

\end{document}